\documentclass[12pt]{article}
\usepackage{color}
\definecolor{darkblue}{rgb}{0.00,0.25,0.50}
\usepackage[colorlinks,filecolor=blue,citecolor=darkblue]{hyperref}

\usepackage[russian, english]{babel}
\usepackage{amssymb,amsmath, amsfonts, mathrsfs, url, amscd}

\begin{document}

\selectlanguage{russian} \setcounter{page}{1} \thispagestyle{empty}

\thispagestyle{empty}

\begin{center}
\textbf{{ПРИБЛИЖЕНИЕ НЕПРЕРЫВНЫХ ПЕРИОДИЧЕСКИХ ФУНКЦИЙ СУММАМИ ВАЛЛЕ ПУССЕНА}}
\end{center}
\vskip0.5cm
\begin{center}
Е.\,Ю.~Овсий, А.\,С. Сердюк
\end{center}
\vskip1cm

\abstract{We obtain an estimate of the deviation of de la Vall\'{e}e Poussin sums
$V_{n,\frac{n}{2}}(f;x)$ from continuous functions $f$, expressed in terms of values of
theirs modulus of continuity. It is established that this estimate can't be improved by
using the well-known analogue of the Lebesgue inequality for  de la Vall\'{e}e Poussin sums.

\vskip0.6cm

Получена оценка уклонений сумм Валле Пуссена $V_{n,\frac{n}{2}}(f;x)$ от непрерывных функций
$f$, выражаемая  через значения их модулей непрерывности. Установлено, что данную оценку
нельзя улучшить путем использования известного аналога неравенства Лебега для сумм Валле
Пуссена.}

\vskip2cm


    Обозначим через $C$ пространство непрерывных $2\pi $-периодических функций $f(x)$, в котором
норма определяется равенством $\|f\|_C=\max\limits_{x}|f(x)|.$

В роботах [\ref{Vallee Poussin}] и [\ref{Vallee Poussin_1919}, с. 33--35] Валле Пуссеном
были рассмотрены суммы
    \begin{equation}\label{2.11.11-16:11:26}
    V_{n,p}(f;x):=\frac{1}{p}\sum\limits_{k=n-p}^{n-1}S_{k}(f;x),\ \ p,n\in \mathbb{N},
    \ \ p\leqslant n,
    \end{equation}
где $S_k(f;x)$ --- частная сумма Фурье  порядка $k$ $2\pi $-периодической функции $f(x).$ Им
же было доказано, что
    \begin{equation}\label{2.11.11-16:12:22}
    |f(x)-V_{n,p}(f;x)|\leqslant 2\frac{n}{p}E_{n-p}(f),
    \end{equation}
где
    $$E_m(f):=\inf_{t_{m-1}}\|f(\cdot)-t_{m-1}(\cdot)\|_C$$

\hskip 2cm

{\footnotesize \noindent $\overline{\ \ \ \ \ \ \ \ \ \ \ \ \ \ \ \ \ \ \ \ \ \ \ \ \ }$\\
\noindent Работа частично поддержана Государственным фондом фундаментальных исследований
Украины (проект $\Phi35/001$)}

\newpage

\hskip -5mm --- наилучшее равномерное приближение функции $f\in C$ тригонометрическими
полиномами $t_{m-1}$, порядок которых не превышает $m-1.$

Впоследствии, суммы вида (\ref{2.11.11-16:11:26}) получили название сумм Валле Пуссена. Как
следует из (\ref{2.11.11-16:12:22}), в случае, когда отношение $\frac{n}{p}$ ограничено,
суммы Валле Пуссена дают наилучший порядок приближения для любой функции $f\in C$,
удовлетворяющей условию \mbox{$E_{n-p}(f)\leqslant KE_n(f).$}

При $p=n$ суммы Валле Пуссена совпадают с суммами Фейера
    $$\sigma _{n-1}(f ;x)=\frac{1}{n}\sum\limits_{k=0}^{n-1}S_k(f ;x),$$
а при $p=1$ --- с суммами Фурье $S_{n-1}(f ;x)$ порядка $n-1.$

Изучению аппроксимативных свойств сумм Валле Пуссена на различных функциональных классах
посвящено значительное количество работ, с результатами которых можно ознакомиться,
например, в работе [\ref{STEANETS-2002/1}].

На классах
    $$H_\omega=\{f\in C:\ \ |f(x')-f(x'')|\leqslant \omega (|x'-x''|),\ \ x',x''\in \mathbb{R}\},$$
где $\omega (\cdot)$ --- заданный модуль непрерывности, точное значение величины
    $$\mathcal{E}(H_\omega ;V_{n,p}):=\sup\limits_{f\in H_\omega }{\|f(\cdot)-V_{n,p}(f;\cdot)\|_C}$$
при $p=n$ получено в [\ref{Stepanets_1981}, с. 222]. Для $p=\overline{1,n-1},$
$p\in\mathbb{N},$ наиболее общие результаты принадлежат \mbox{А.В. Ефимову.} В частности, в
\mbox{[\ref{Efimov_1959}, с. 739]} им доказано следующее соотношение:
\begin{equation}\label{3.11.11-12:01:58}
    \mathcal{E}(H_\omega ;V_{n,p})=A_{n,p}(\omega )+O(1)\omega (1/n),
\end{equation}
где
    $$A_{n,p}(\omega )=
  \begin{cases}
    \frac{e_n(\omega )}{\pi }\ln\frac{n-1}{p}, & 1\leqslant p\leqslant \frac{n+1}{2}, \\
    \frac{2}{\pi p}\int\limits_{\frac{1}{n}}^{\frac{1}{n-p}}\frac{\omega (t)}{t^2}\,dt,
    & \frac{n+1}{2}\leqslant p\leqslant n-1,
  \end{cases}
    $$
а $e_n(\omega )$ --- верхняя грань по классу $H_\omega $ $n$-го коэффициента Фурье.
 Из (\ref{3.11.11-12:01:58})
следует, что при $p=\frac{n}{2}$
    \begin{equation}\label{3.11.11-12:30:05}
    \mathcal{E}(H_\omega ;V_{n,\frac{n}{2}})\leqslant K\omega (1/n),
    \end{equation}
то есть суммы Валле Пуссена реализуют порядок наилучшего приближения на $H_\omega .$ О
величине константы $K$ в (\ref{3.11.11-12:30:05}) можно судить на основании аналога
неравенства Лебега для сумм Валле Пуссена (см., например, [\ref{Stechkin_78}, с. 61])
\begin{equation}\label{3.11.11-12:34:52}
    |f(x)-V_{n,p}(f;x)|\leqslant
    \Big(\sup\limits_{\|f\|_C\leqslant 1}{\|V_{n,p}(f;\cdot)\|_C}+1\Big)E_{n-p}(f).
\end{equation}
Как вытекает из работы [\ref{Stechkin_51}],
    $$\sup\limits_{\|f\|_C\leqslant 1}{\|V_{n,\frac{n}{2}}(f;\cdot)\|_C}=\frac{
    1}{2\pi }\int\limits_{-\pi }^{\pi }\bigg|\frac{\sin\frac{3}{2}t}{\sin\frac{t}{2}}\bigg|\,dt
    =\frac{1}{3}+\frac{2\sqrt{3}}{\pi }.$$
Кроме того (см., например, [\ref{KORNEICHUK_Exact_Constant}, с. 261]), для любой функции
$f\in C$ \mbox{$(f\neq \text{const})$} выполняется неравенство
\begin{equation}\label{3.11.11-12:59:10}
    E_m(f)<\omega \Big(f,\frac{\pi }{m}\Big),\ \ \ m=1,2,\ldots,
\end{equation}
а для любой функции $f\in C_*$, где $C_*$ --- подмножество функций $f$ из $C,$ у которых
модуль непрерывности $\omega (f,t)$ является выпуклым вверх, --- неравенство
\begin{equation}\label{3.11.11-13:00:53}
    E_m(f)\leqslant \frac{1}{2}\,\omega \Big(f,\frac{\pi }{m}\Big),\ \ \ m=1,2,\ldots.
\end{equation}
В силу (\ref{3.11.11-12:34:52}) и (\ref{3.11.11-12:59:10}) для любой функции $f\in C$ при
 $p=\frac{n}{2}$ получаем
\begin{equation}\label{10.11.11-22:29:14}
    |f(x)-V_{n,\frac{n}{2}}(f;x)|<\Big(\frac{4}{3}+\frac{2\sqrt{3}}{\pi }\Big)\,\omega
    \Big(f,\frac{2\pi }{n}\Big),
\end{equation}
а в силу (\ref{3.11.11-12:34:52}) и (\ref{3.11.11-13:00:53}) для любой функции $f\in C_*$
будем иметь
\begin{equation}\label{10.11.11-22:33:40}
    |f(x)-V_{n,\frac{n}{2}}(f;x)|\leqslant \Big(\frac{2}{3}+\frac{\sqrt{3}}{\pi }\Big)\omega
    \Big(f,\frac{2\pi }{n}\Big).
\end{equation}
Если $f\in H_\omega $, то на основании (\ref{10.11.11-22:29:14}) и (\ref{10.11.11-22:33:40})
заключаем, что
\begin{equation}\label{3.11.11-13:06:55}
    \mathcal{E}(H_\omega ;V_{n,\frac{n}{2}})<\Big(\frac{4}{3}+\frac{2\sqrt{3}}{\pi }\Big)\omega
    \Big(\frac{2\pi }{n}\Big),
\end{equation}
если $\omega(\cdot) $ --- произвольный модуль непрерывности, и
\begin{equation}\label{3.11.11-13:17:03}
    \mathcal{E}(H_\omega ;V_{n,\frac{n}{2}})\leqslant \Big(\frac{2}{3}+\frac{\sqrt{3}}{\pi }\Big)\omega
    \Big(\frac{2\pi }{n}\Big),
\end{equation}
если $\omega(\cdot) $ --- выпуклый модуль непрерывности. Возникает вопрос: являются ли
оценки (\ref{10.11.11-22:29:14}) и (\ref{10.11.11-22:33:40}), а соответственно
(\ref{3.11.11-13:06:55}) и (\ref{3.11.11-13:17:03}), неулучшаемыми? В дальнейшем будет
показано, что в отношении (\ref{10.11.11-22:29:14}), (\ref{3.11.11-13:06:55}) и
(\ref{3.11.11-13:17:03}) это не так.

Имеет место утверждение.

\textbf{Теорема 1.} \emph{Пусть $f\in C$. Тогда}
\begin{equation}\label{10.11.11-22:43:12}
    |f(x)-V_{n,\frac{n}{2}}(f;x)|<
    \omega \Big(f,\frac{6\pi }{7n}\Big)+
    \frac{9}{10\pi }\,\omega \Big(f,\frac{2\pi }{3n}\Big)+
    \frac{31}{25\pi }\,\omega \Big(f,\frac{\pi }{n}\Big).
\end{equation}

Поскольку
    \begin{equation}\label{10.11.11-22:54:21}
    \omega \Big(f,\frac{6\pi }{7n}\Big)+
    \frac{9}{10\pi }\,\omega \Big(f,\frac{2\pi }{3n}\Big)+
    \frac{31}{25\pi }\,\omega \Big(f,\frac{\pi }{n}\Big)<1.6812\,\omega
    \Big(f,\frac{\pi }{n}\Big)<$$
    $$<
    \Big(\frac{4}{3}+\frac{2\sqrt{3}}{\pi }\Big)\omega \Big(f,\frac{\pi }{n}\Big),
    \end{equation}
то очевидно, что оценка (\ref{10.11.11-22:43:12}) уточняет (\ref{10.11.11-22:29:14}).

Если $f\in H_\omega $,  то из (\ref{10.11.11-22:43:12}) получаем утверждение.

\textbf{\emph{Следствие} 1.} \emph{Пусть $\omega (\cdot)$ --- произвольный модуль
непрерывности. Тогда}
\begin{equation}\label{6.11.11-12:51:45}
    \mathcal{E}(H_\omega ;V_{n,\frac{n}{2}})<\omega \Big(\frac{6\pi }{7n}\Big)+
    \frac{9}{10\pi }\,\omega \Big(\frac{2\pi }{3n}\Big)+\frac{31}{25\pi }\,\omega \Big(\frac{\pi }{n}\Big).
\end{equation}

Учитывая (\ref{10.11.11-22:54:21}), приходим к выводу, что (\ref{6.11.11-12:51:45}) является
более точной оценкой величины $\mathcal{E}(H_\omega ;V_{n,\frac{n}{2}})$, чем
(\ref{3.11.11-13:06:55}). Заметим также, что в случае, когда $H_\omega $ является классом
Гельдера $H^\alpha $ (т.е. $\omega (t)=t^\alpha $), оценка (\ref{6.11.11-12:51:45}) уточняет
(\ref{3.11.11-13:17:03}) при всех $\alpha \in[0.38,1].$

Для класса $H^\alpha $ можно получить двухстороннюю оценку величины $\mathcal{E}(H^\alpha
;V_{n,\frac{n}{2}})$. Для этого воспользуемся установленной в работе
\mbox{[\ref{KORNEICHUK_1963}, с. 42]}  оценкой
    \begin{equation}\label{6.11.11-14:04:28}
    \sup\limits_{f\in H^\alpha  }\|f(\cdot)-U_{n-1}(f;\cdot)\|_C\geqslant
    \frac{\pi ^\alpha }{(1+\alpha )n^\alpha },\ \ n\in \mathbb{N},
    \end{equation}
в которой
\begin{equation}\label{9.11.11-12:50:11}
U_{n-1}(f ;x)=\frac{a_0(f)}{2}+\sum\limits_{ k=1}^{n-1}\lambda _{k,n}(a_k(f )\cos kx+b_k(f
)\sin kx),
\end{equation}
$\lambda _{k,n}$ --- произвольная заданная последовательность действительных чисел, а $a_k(f
)$ и $b_k(f )$ --- коэффициенты Фурье функции $f(x) $. Поскольку суммы Валле Пуссена
$V_{n,p}(f ;x)$ допускают представление в виде (\ref{9.11.11-12:50:11}), то из теоремы 1 и
оценки (\ref{6.11.11-14:04:28}) получаем следующее утверждение.

\textbf{Утверждение 1.} \emph{При $0<\alpha \leqslant 1$ имеет место двухсторонняя оценка}
    \begin{equation}\label{9.11.11-23:23:21}
    \frac{1}{(\alpha +1)}\frac{\pi ^\alpha }{n^\alpha }\leqslant \mathcal{E}(H^\alpha  ;V_{n,\frac{n}{2}}
    )<\bigg(\frac{6^\alpha }{7^\alpha }+\frac{3^{2-\alpha }2^{\alpha -1}}{5\pi } +
    \frac{31}{25\pi }\bigg)\frac{\pi ^\alpha }{n^\alpha }.
    \end{equation}

В частности, при $\alpha =1$ из (\ref{9.11.11-23:23:21}) получаем
    $$\frac{\pi }{2n}\leqslant \mathcal{E}(H^1  ;V_{n,\frac{n}{2}}
    )<\frac{1.443\pi }{n},$$
а при $\alpha =\frac{1}{2}$
    $$\frac{2}{3}\sqrt{\frac{\pi }{n}}\leqslant \mathcal{E}(H^{\frac{1}{2}}  ;V_{n,\frac{n}{2}})<
    1.555\sqrt{\frac{\pi }{n}}.$$

\textbf{\emph{Доказательство теоремы} 1.} Рассмотрим уклонение
\begin{equation}\label{6.10.11-14:51:00}
    \rho _{n,p}(f ;x):=f (x)-V_{n,p}(f ;x),\ \ x\in \mathbb{R}.
\end{equation}
Поскольку для произвольной функции $f \in C$ и $p=\overline{1,n}$ имеет место представление
\begin{equation}\label{6.10.11-15:46:43}
    \rho _{n,p}(f ;x)=\frac{n}{\pi p}\int\limits_{-\infty}^{\infty}
    \bigg(f (x)-f \Big(x+\frac{t}{n}\Big)\bigg)
    \frac{\cos\frac{n-p}{n}t-\cos t}{t^2}\,dt,\ \ x\in\mathbb{R},
\end{equation}
(см., например, [\ref{Efimov_1959}, с. 750]), то, полагая в (\ref{6.10.11-15:46:43})
$p=\frac{n}{2}$, получаем
\begin{equation}\label{6.10.11-15:58:02}
    \rho _{n,\frac{n}{2}}(f ;x)=\frac{2}{\pi}\int\limits_{-\infty}^{\infty}
    \bigg(f (x)-f \Big(x+\frac{t}{n}\Big)\bigg)
    \frac{\cos\frac{t}{2}-\cos t}{t^2}\,dt=$$
    $$=\frac{2}{\pi}\int\limits_{0}^{\infty}
    \bigg(2f (x)-f \Big(x+\frac{t}{n}\Big)-f \Big(x-\frac{t}{n}\Big)\bigg)
    \frac{\cos\frac{t}{2}-\cos t}{t^2}\,dt.
\end{equation}

Покажем, что на каждом интервале $(x_k,x_{k+1}),$ $k\in\mathbb{N}$, $x_k:=$\mbox{$(2k-1)\pi
$,} функция
\begin{equation}\label{6.10.11-16:01:29}
    g(x):=2\int\limits_{x}^{\infty}\frac{\cos\frac{t}{2}-\cos t}{t^2}\,dt, \ \ x>0,
\end{equation}
имеет, по крайней мере, один ноль. Интегрируя по частям, находим
\begin{equation}\label{6.10.11-16:08:54}
    g(x)=\frac{1}{2x}\Big(
    4\cos\frac{x}{2}+2x\text{Si}\Big(\frac{x}{2}\Big)-4\cos x-4x\text{Si}(x)+\pi x\Big),
\end{equation}
где $\text{Si}(x)$ --- интегральный синус, то есть функция вида
    $$\text{Si}(x)=\int\limits_{0}^{x}\frac{\sin t}{t}\,dt.$$
С учетом равенства
    $$\text{Si}(x)=\frac{\pi }{2}-\frac{\cos x}{x}+\int\limits_{
    x}^{\infty}\frac{\cos t}{t^2}\,dt,$$
 из (\ref{6.10.11-16:08:54}) получаем
\begin{equation}\label{6.10.11-16:13:17}
    g(x)=2\int\limits_{x/2}^{x}\frac{\cos t}{t^2}\,dt-\int\limits_{
    x/2}^{\infty}\frac{\cos t}{t^2}\,dt.
\end{equation}
Покажем, что
\begin{equation}\label{6.10.11-16:16:49}
    \bigg|\int\limits_{
    x_k/2}^{\infty}\frac{\cos t}{t^2}\,dt\bigg|<2\bigg|\int\limits_{x_k/2}^{x_k}\frac{\cos t}{t^2}\,dt\bigg|,\ \ x_k=(2k-1)\pi ,\ \ k\in\mathbb{N}.
\end{equation}
Легко видеть, что
\begin{equation}\label{6.10.11-16:20:42}
    \bigg|\int\limits_{
    x_k/2}^{\infty}\frac{\cos t}{t^2}\,dt\bigg|<\bigg|\int\limits_{
    x_k/2}^{x_k/2+\pi }\frac{\cos t}{t^2}\,dt\bigg|.
\end{equation}
Если $k=1$, то очевидно, что
\begin{equation}\label{6.10.11-16:22:41}
    \bigg|\int\limits_{
    x_k/2}^{x_k/2+\pi }\frac{\cos t}{t^2}\,dt\bigg|<2\bigg|\int\limits_{x_k/2}^{x_k}\frac{\cos t}{t^2}\,dt\bigg|
\end{equation}
и отсюда, с учетом (\ref{6.10.11-16:20:42}), сразу следует (\ref{6.10.11-16:16:49}).

Пусть $k\geqslant 2.$ Поскольку
    $$\text{sign}\bigg(
    \int\limits_{x_k/2}^{(x_k+\pi )/2}\frac{\cos t}{t^2}\,dt\bigg)=(-1)^k,\ \ k=2,3,\ldots,$$
и
    $$\text{sign}\bigg(\int\limits_{
    (x_k+\pi )/2}^{x_k}\frac{\cos t}{t^2}\,dt\bigg)=\text{sign}
    \bigg(\int\limits_{k\pi }^{(2k-1)\pi }
    \frac{\cos t}{t^2}\,dt\bigg)=$$
    $$=\text{sign}\bigg(
    \int\limits_{k\pi }^{(2k-1)\pi }\frac{\sin t}{t^3}\,dt\bigg)
    =(-1)^k,\ \ k=2,3,\ldots,$$
то
\begin{equation}\label{6.10.11-16:36:33}
    2\bigg|\int\limits_{x_k/2}^{x_k}\frac{\cos t}{t^2}\,dt\bigg|>
    2\bigg|\int\limits_{x_k/2}^{(x_k+\pi )/2}\frac{\cos t}{t^2}\,dt\bigg|>
    \bigg|\int\limits_{x_k/2}^{x_k/2+\pi }\frac{\cos t}{t^2}\,dt\bigg|,\ \ k=2,3,\ldots
\end{equation}
Объединяя (\ref{6.10.11-16:20:42}) и (\ref{6.10.11-16:36:33}) приходим к
(\ref{6.10.11-16:16:49}), из которого следует, что
\begin{equation}\label{6.10.11-16:57:32}
    \text{sign}(g(x_k))=\text{sign}\bigg(
    2\int\limits_{x_k/2}^{x_k}\frac{\cos t}{t^2}\,dt-\int\limits_{
    x_k/2}^{\infty}\frac{\cos t}{t^2}\,dt\bigg)=$$
    $$=\text{sign}\bigg(
    \int\limits_{x_k/2}^{x_k}\frac{\cos t}{t^2}\,dt\bigg)=(-1)^k,\ \ x_k=(2k-1)\pi ,\ \ k\in\mathbb{N}.
\end{equation}

Обозначим через $\tau _k$ ноль функции $g(x)$ на интервале $(x_{k-1},x_{k}),$
\mbox{$x_k=(2k-1)\pi,$} $k\in\mathbb{N},$ $x_0=0.$ Тогда, используя обозначения
    $$\tau _0=0,$$
    $$\Delta (f;x;t/n):=2f (x)-f \Big(x+\frac{t}{n}\Big)-f
    \Big(x-\frac{t}{n}\Big),$$
из (\ref{6.10.11-15:58:02}) получаем
    $$|\rho _{n,\frac{n}{2}}(f ;x)|=\frac{
    2}{\pi }\bigg|
    \sum\limits_{k=0}^{\infty}\int\limits_{\tau _k}^{\tau _{k+1}}\Delta (f;x;t/n)
    \frac{\cos\frac{t}{2}-\cos t}{t^2}\,dt\bigg|.$$
Далее, поскольку
    $$
    \int\limits_{\tau _k}^{\tau _{k+1}}\frac{\cos\frac{t}{2}-\cos t}{t^2}\,dt=
    \frac{g(\tau _k)-g(\tau _{k+1})}{2}=0,$$
приходим к неравенству
\begin{equation}\label{8.11.11-19:05:53}
    |\rho _{n,\frac{n}{2}}(f ;x)|\leqslant \frac{2}{\pi }
    \bigg(\bigg|\int\limits_{\tau _0}^{\tau _1}\Delta (f;x;t/n)\frac{\cos \frac{t}{2}-\cos t}{t^2}\,dt\bigg|+$$
    $$+
    \bigg|\sum\limits_{k=1}^{\infty}\int\limits_{\tau _k}^{\tau _{k+1}}\Big(\Delta (f;x;t/n)
    -\Delta (f;x;\frac{\tau _k+\tau _{k+1}}{2n})\Big)
    \frac{\cos\frac{t}{2}-\cos t}{t^2}\,dt\bigg|
    \bigg).
\end{equation}
Учитывая, что для любой функции $f\in C $ выполняются неравенства:
    $$|\Delta (f;x;t/n)|\leqslant 2\omega (f,t/n),\ \ t>0$$
и
    $$|\Delta (f;x;t_2/n)-\Delta (f;x;t_1/n)|
    \leqslant 2\omega \Big(f,\frac{|t_2-t_1|}{n}\Big),\ \ t_1,t_2>0,$$
из (\ref{8.11.11-19:05:53}) получаем
    $$|\rho _{n,\frac{n}{2}}(f ;x)|\leqslant \frac{4}{\pi }
    \bigg(\int\limits_{\tau _0}^{\tau _1}\omega
    (f,t/n)\frac{|\cos \frac{t}{2}-\cos t|}{t^2}\,dt+$$
    \begin{equation}\label{10.10.11-16:30:35}
    +
    \sum\limits_{k=1}^{\infty}\int\limits_{
    \tau _k}^{\tau _{k+1}}\omega \Big(f,
    \frac{|2t-\tau _k-\tau _{k+1}|}{2n}\Big)\frac{|\cos \frac{t}{2}-\cos t|}{t^2}\,dt\bigg).
    \end{equation}
Вычисления, произведенные на ЭВМ, показывают, что для первых пяти нулей функции $g(x)$
справедливы следующие включения:
    \begin{equation}\label{8.11.11-18:57:31}
    \tau _1\in(2.657, 2.66),
    \end{equation}
    \begin{equation}\label{8.11.11-18:57:32}
    \tau _2\in(6.83, 6.84),
    \end{equation}
    \begin{equation}\label{8.11.11-18:57:33}
    \tau _3\in(14.16, 14.17),
    \end{equation}
    \begin{equation}\label{8.11.11-18:57:35}
    \tau _4\in(19.09, 19.10),
    \end{equation}
    \begin{equation}\label{8.11.11-18:57:36}
    \tau _5\in(26.41, 26.42).
    \end{equation}
Тогда, на основании (\ref{8.11.11-18:57:31}) для первого интеграла в правой части
(\ref{10.10.11-16:30:35}) можем записать
\begin{equation}\label{8.11.11-19:15:09}
    \int\limits_{\tau _0}^{\tau _1}\omega
    (f,t/n)\frac{|\cos \frac{t}{2}-\cos t|}{t^2}\,dt\leqslant \omega (f,\tau _1/n)
    \int\limits_{\tau _0}^{\tau _1}\frac{|\cos \frac{t}{2}-\cos t|}{t^2}\,dt<$$
    $$<\omega \Big(f,\frac{2.66}{n}\Big)\int\limits_{0}^{2.66}\frac{|\cos \frac{t}{2}-
    \cos t|}{t^2}\,dt<0.786\,\omega \Big(f,\frac{2.66}{n}\Big)<$$
    $$<0.786\,\omega \Big(f,\frac{6\pi }{7n}\Big).
\end{equation}
Аналогичным образом, учитывая (\ref{8.11.11-18:57:32})--(\ref{8.11.11-18:57:36}), находим
    $$\sum\limits_{k=1}^{4}\int\limits_{
    \tau _k}^{\tau _{k+1}}\omega \Big(f,
    \frac{|2t-\tau _k-\tau _{k+1}|}{2n}\Big)\frac{|\cos \frac{t}{2}-\cos t|}{t^2}\,dt\leqslant $$
    $$\leqslant
    \sum\limits_{k=1}^{4}\omega \Big(f,\frac{\tau _{k+1}-\tau _k}{2n}\Big)
    \int\limits_{\tau _k}^{\tau _{k+1}}\frac{|\cos\frac{t}{2}-\cos t|}{t^2}\,dt<$$
    \begin{equation}\label{6.11.11-15:06:32}
    <0.225\,\omega \Big(f,\frac{2.0915}{n}\Big)+0.057\,\omega \Big(f,\frac{3.67}{n}\Big)+
    0.019\,\omega \Big(f,\frac{2.47}{n}\Big)+$$$$
    +0.011\,\omega \Big(f,\frac{3.67}{n}\Big)<0.225\,\omega \Big(f,\frac{2\pi }{3n}\Big)+0.068\,\omega
    \Big(f,\frac{3.67}{n}\Big)+0.019\,\omega \Big(f,\frac{\pi }{n}\Big).
    \end{equation}
Поскольку
    $$\omega
    \Big(f,\frac{3.67}{n}\Big)\leqslant \omega
    \Big(f,\frac{3.67-\pi }{n}\Big)+\omega\Big(f,\frac{\pi }{n}\Big)<2\omega
    \Big(f,\frac{\pi }{n}\Big),$$
то из (\ref{6.11.11-15:06:32}) получаем
\begin{equation}\label{8.11.11-19:41:10}
    \sum\limits_{k=1}^{4}\int\limits_{
    \tau _k}^{\tau _{k+1}}\omega \Big(f,
    \frac{|2t-\tau _k-\tau _{k+1}|}{2n}\Big)\frac{|\cos \frac{t}{2}-\cos t|}{t^2}\,dt<$$
    $$<
    0.225\,\omega \Big(f,\frac{2\pi }{3n}\Big)+0.155\,\omega \Big(f,\frac{\pi }{n}\Big).
\end{equation}
Далее, в силу того, что
    $$\tau _{k+1}-\tau _k<x_{k+1}-x_{k-1}=4\pi ,\ \ \ x_k=(2k-1)\pi,\ \ \ k\in\mathbb{N},$$
имеем
\begin{equation}\label{10.10.11-16:34:41}
    \sum\limits_{k=5}^{\infty}\int\limits_{
    \tau _k}^{\tau _{k+1}}\omega \Big(f,
    \frac{|2t-\tau _k-\tau _{k+1}|}{2n}\Big)\frac{|\cos \frac{t}{2}-\cos t|}{t^2}\,dt\leqslant
    $$
    $$\leqslant 2\sum\limits_{k=5}^{\infty}\omega \Big(f,\frac{\tau _{k+1}-\tau _k}{2n}\Big)\int\limits_{\tau _k}^{\tau _{k+1}}
    \frac{dt}{t^2}<2\omega \Big(f,\frac{2\pi }{n}
    \Big)\int\limits_{\tau _5}^{\infty}\frac{dt}{t^2}<0.152\,\omega \Big(f,\frac{\pi }{n}\Big).
\end{equation}
Объединяя (\ref{10.10.11-16:30:35}), (\ref{8.11.11-19:15:09}), (\ref{8.11.11-19:41:10}) и
(\ref{10.10.11-16:34:41}), приходим к оценке
\begin{equation}\label{10.10.11-16:48:29}
    |\rho _{n,\frac{n}{2}}(f ;x)|<\frac{4}{\pi }\bigg(
    0.768\,\omega \Big(f,\frac{6\pi }{7n}\Big)+0.225\,\omega \Big(f,\frac{2\pi }{3n}\Big)+0.307
    \,\omega \Big(f,\frac{\pi }{n}\Big)\bigg)<$$
    $$<\omega \Big(f,\frac{6\pi }{7n}\Big)+\frac{9}{10\pi }\,
    \omega \Big(f,\frac{2\pi }{3n}\Big)+\frac{31}{25\pi }\,\omega \Big(f,\frac{\pi }{n}\Big),\ \ x\in \mathbb{R},
\end{equation}
откуда непосредственно следует справедливость теоремы 1.

\vskip 3.5mm

\footnotesize
\begin{enumerate}
\Rus

\item\label{Vallee Poussin}{\it{De La Vall\'{e}e Poussin C.}\/} Sur la meilleure
approximation des fonctions d'une variable r\'{e}elle par des expressions d'ordre donn\'{e}
// Comptes
rendus de l'Acad\'{e}mie des Sciences. --- 1918. --- \textbf{166}, \mbox{№ 4.} --- P.
799---802.

\item\label{Vallee Poussin_1919}{\it{De La Vall\'{e}e Poussin C.}\/} Le\c{c}ons sur
l'approximation des fonctions d'une variable r\'{e}elle. --- Paris: Gauthier-Villars, 1919.
--- 152 p.

\item\label{STEANETS-2002/1}{\it{Степанец А.И.}\/} Методы теории приближений: В 2 ч. //
Праці Ін-ту математики НАН України. --- 2002. --- \textbf{{40}}. --- Ч. IІ. --- 468 с.

\item\label{Stepanets_1981}{\it{Степанец А.И.}\/} Равномерные приближения
тригонометрическими полиномами. --- Киев: Наук. думка, 1981. --- 340 с.

\item\label{Efimov_1959}{\it{Ефимов А.В.}\/} О приближении периодических функций суммами
Валле Пуссена
// Изв. АН СССР. --- 1959. --- \textbf{{23}}. --- С.
737---770.

\item\label{Stechkin_78}{\it{Ste\v{c}kin S.B.}\/} On the approximation of periodic functions
by de la Vall\'{e}e Poussin sums // Anal. Math. --- 1978. --- \textbf{4}. --- С. 61---74.

\item\label{Stechkin_51}{\it{Стечкин С.Б.}\/} О суммах Валле-Пуссена // Докл. АН СССР.
--- 1951. --- \textbf{80}.--- № 4. --- С. 545---548.

\item\label{KORNEICHUK_Exact_Constant}{\it{Корнейчук Н.П.}\/} Точные константы в теории
приближения. --- М.: Наука, 1987. --- 424 с.

\item\label{KORNEICHUK_1963}{\it{Корнейчук Н.П.}\/} О наилучшем приближении непрерывных
функций // Изв. АН СССР. --- 1963. --- \textbf{27}. --- С. 29---44.

\end{enumerate}

\label{end}

\emph{Contact information}: \href{http://www.imath.kiev.ua/~funct}{Department of the Theory
of Functions}, Institute of Mathematics of Ukrainian National Academy of Sciences, 3,
Tereshenkivska st., 01601, Kyiv, Ukraine \vskip 0.2 cm \emph{E-mail}:
\href{mailto:ievgen.ovsii@gmail.com}{ievgen.ovsii@gmail.com},
\href{mailto:serdyuk@imath.kiev.ua}{serdyuk@imath.kiev.ua}

\end{document}